\documentclass[12pt]{amsart}
\usepackage[active]{srcltx}
\def\qed{\hfill \rule{4pt}{7pt}}
\def\pf{\noindent {\it Proof.} }
\newcommand{\poq}[2]{(#1;q)_{#2}}
\newcommand{\poqt}[2]{(#1;q^{-1})_{#2}}
\def\qed{\hfill \rule{4pt}{7pt}}
\def\pf{\noindent {\it Proof.} }

\usepackage{latexsym}
\usepackage{amsmath,amsthm}
\usepackage{graphicx}
\usepackage{amssymb,amsmath}

\theoremstyle{definition}
\newtheorem{dl}{Theorem}
\newtheorem{tl}{Corollary}
\newtheorem{yl}{Lemma}

\begin{document}

\title[A  Five-variable  Ramanujan's reciprocity theorem]{A  Five-variable generalization of\\
 Ramanujan's reciprocity theorem\\
  and its applications}
\author{X. R. Ma}
\subjclass[2000]{Primary 05A10,33D15}
\thanks{ Supported by  NSFC  grant No. 10771156.}
\address{The Corresponding Author       \newline
    Dr. X.Ma         \newline
    Department of Matematics\newline
        SuZhou University\newline
        SuZhou\: \:P.\:O.\:Box\:173 \newline
        215006 SuZhou,$\quad$P.R.China       \newline
       Email \emph{xrma@public1.sz.js.cn}}
\maketitle \vskip 5mm
\date{\today}
\maketitle
\begin{abstract}
    By virtue of  Bailey's  well-known bilateral
$\,_6\psi_6$ summation formula and Watson's transformation formula,
we extend the four-variable generalization of Ramanujan's
reciprocity theorem due to Andrews to a five-variable one. Some
relevant new $q$-series identities  including a new proof of
Ramanujan's reciprocity theorem and
of Watson's  quintuple product identity only based  on  Jackson's transformation  are presented. \vspace{10pt}\\
Keywords: $q$-series; reciprocity theorem; Ramanujan's $\,_1\psi_1$
summation formula; Jacobi's triple product identity; Watson's
quintuple product identity; transformation formula.
\end{abstract}
\parskip 7pt
\section{Introduction}

In his lost notebook \cite[p.40]{raman}, Ramanujan offered without
proof a beautiful $q$-series identity, which is now called
\textit{Ramanujan's reciprocity theorem}.
\begin{dl}\label{1.1} For  $ a, b \neq q^{-n}$,
 it
holds
\begin{eqnarray}
&&\rho(a,b)-\rho(b,a)=\left(\frac{1}{b}-\frac{1}{a}\right)\frac{\poq{q,aq/b,bq/a}{\infty}}
{\poq{-aq,-bq}{\infty}},
  \end{eqnarray}
  where
\begin{eqnarray}\rho(a,b)=\left(1+\frac{1}{b}\right)\sum_{k=0}^{\infty}
\frac{(-1)^kq^{k(k+1)/2}} {\poq{-aq}{k}} \left(\frac{a}{b}\right)^k.
\end{eqnarray}
\end{dl}

Ramanujan's reciprocity theorem  has been proved to be very useful
to partial theta function identities.
 For further details on this subject, the reader can  refer
the forthcoming second volume of Ramanujan's lost notebook
\cite{ber0} by Andrews and Berndt. It has been an active topic, in
the past years, to find possibly short and easy proofs for Theorem
1.1. Up to now, various approaches have been found by many
mathematicians. For our purpose of this paper, we only mention  a
few remarkable results. As is known to us, Andrews \cite[Theorem
1]{andrews} gave the first proof of Theorem \ref{1.1}, whereas it
seems a bit complicated, by setting a four-variable generalization
of Theorem \ref{1.1} \cite[Theorem 6]{andrews}.   In 2003, Liu
\cite[Theorem 6]{liu} showed the four-variable generalization of
Andrews by  $q$-exponential operator identity.  During the last two
years, Adiga-Anitha \cite{ad}, Berndt et al \cite{ber} found
independently that this theorem can be verified by Heine's
transformation for $\,_2\phi_1$ series. Besides, Berndt et al
presented another simpler analytic proof and an elegant
combinatorial proof in \cite{ber}. An interesting phenomenon   is
that almost all known analytic proofs utilized Ramanujan's
$\,_1\psi_1$
 summation formula   and the Rogers-Fine identity. In a very recent paper  \cite{kang},
  Kang rederived  Andrews' four-variable
generalization from Sears' three-term relation between
$\,_3\phi_2$-series \cite[III.33]{10}, which had been claimed as a
long-waited problem proposed by Andrews and Agarwal.

  We now restate the four-variable
generalization of Theorem \ref{1.1} in Kang's form for  further
discussion.
\begin{dl}\label{1.2} For four parameters $a,b,c,d$ satisfying $c,d\neq
-aq^{-m},-bq^{-n},$ $n,m\geq 0$, $0<|d|<|b|$, it holds
 \begin{align}\label{kang1}
&&\rho(a,b;c,d)-\rho(b,a;c,d)\\
&&=\left(\frac{1}{b}-\frac{1}{a}\right)\frac{\poq{q,aq/b,bq/a,c,d,cd/(ab)
}{\infty}} {\poq{-aq,-bq,-c/b,-d/b,-c/a,-d/a}{\infty}},\nonumber
   \end{align}
  where
\begin{eqnarray}&&\rho(a,b;c,d)=\left(1+\frac{1}{b}\right)\sum_{k=0}^{\infty} \frac{\poq{c,-aq/d}{k}}
 {\poq{-aq}{k}\poq{-c/b}{k+1}}
\left(-\frac{d}{b}\right)^k\label{3a}\\
&=&\left(1+\frac{1}{b}\right)\sum_{k=0}^{\infty}q^{\binom{k+1}{2}}
 \frac{(1+cdq^{2k}/b)\poq{c,d,cd/(ab)}{k}}
{\poq{-aq}{k}\poq{-c/b,-d/b}{k+1}} \left(-\frac{a}{b}\right)^k.
\label{3c}
 \end{eqnarray}
\end{dl}

For the connections of this four-variable reciprocity theorem with
some well known $q$-series identities such as  Jacobi's triple
 and Watson's quintuple product identity, the reader can consult
\cite{kang}.  Here, we  only point out that Identity (\ref{kang1})
given by Kang is  equivalent to  the four-variable generalization of
Theorem \ref{1.1}  of Andrews \cite[Theorem 6]{andrews}. This fact
is stated clearly by \,(4.10) of \cite{kang}. One of the most
important results of Kang, in the author's point of view, is that
she established the new expression \,(\ref{3c}) for the function
$\rho(a,b;c,d).$

In the present paper, motivated by the  method  of Kang, we will
extend Theorem \ref{1.2} to the following five-variable form.

\begin{dl}\label{1.3} For five parameters $a,b,c,d,e$ satisfying
\[0<|cde|<|abq|,c,d,e\neq
-aq^{-m},-bq^{-n},n,m\geq 0,\] it holds
 \begin{align}
&&\rho(a,b;c,d,e)-\rho(b,a;c,d,e)=\left(\frac{1}{b}-\frac{1}{a}\right)\\
&&\times\frac{\poq{q,aq/b,bq/a,c,d,e,cd/(ab),ce/(ab), de/(ab)
}{\infty}} {\poq{-aq,-bq,-c/a,-c/b,-d/a,-d/b,-e/a,-e/b, cde/(abq)
}{\infty}},\nonumber
   \end{align}
  where
\begin{eqnarray}\rho(a,b;c,d,e)&=&\sum_{k=0}^{\infty}
\left(1-\frac{aq^{2k+1}}{b}\right)\frac{\poq{-1/b}{k+1}}
{\poq{-c/b,-d/b,-e/b}{k+1}}\nonumber\\
&\times&\frac{\poq{-aq/c,-aq/d,-aq/e}{k}} {\poq{-aq}{k}}
\left(\frac{cde}{ab q}\right)^k.
 \end{eqnarray}
\end{dl}
As application of Theorem \ref{1.3}, several new $q$-series
identities will also be derived.

To make our paper self-contained,  we will repeat a few  standard
notation and terminology for basic hypergeometric series (or
$q$-series) found in \cite{10}. Given a (fixed) complex number $q$
with $|q|<1$, a complex number $a$ and an integer $n$, define the
$q$-shifted factorials
 $(a;q)_\infty$ and $(a; q)_n$  as
\begin{eqnarray}
(a;q)_{\infty}=\prod_{k=0}^{\infty}(1-aq^k),\quad (a;q)_n=
\frac{\poq{a}{\infty}}{\poq{aq^n}{\infty}}. \label{conven}
    \end{eqnarray}
We also employ the following compact multi-parameter notation
\begin{eqnarray*}
  (a_1,a_2,\cdots,a_m;q)_n = (a_1;q)_n(a_2;q)_n\cdots
   (a_m;q)_n.
    \end{eqnarray*}  The basic and bilateral hypergeometric series with the
base $q$  are defined respectively as
\begin{eqnarray*}&&{}_{r}\phi _{r-1}\left[\begin{matrix}a_{1},\dots ,a_{r}
\\ b_{1},\dots ,b_{r-1}\end{matrix}
; q, z\right]=\sum _{n=0} ^{\infty }\frac{\poq {a_{1},\cdots
,a_{r}}{n}}{\poq
{q,b_{1},\cdots,b_{r-1}}{n}}z^{n};\\
&&{}_{r}\psi _{r}\left[\begin{matrix}a_{1},\dots ,a_{r}\\
b_{1},\dots ,b_{r}\end{matrix} ; q, z\right]=\sum _{n=-\infty}
^{\infty }\frac{\poq {a_{1},\cdots ,a_{r}}{n}}{\poq
{b_{1},\cdots,b_{r}}{n}}z^{n}.
\end{eqnarray*}
In addition, the compact notation
$\,_{r}W_{r-1}(a_1;a_4,\cdots,a_{r};q,z)$ denotes the special case
 of the above ${}_{r}\phi _{r-1}$ called  a very-well-poised series, in which all parameters
satisfy the relations
$$
qa_1=b_1a_2=\cdots=b_{r-1}a_{r}; a_2=q\sqrt{a_1},a_3=-q\sqrt{a_1}.
$$

\section{Proof of the main result}\setcounter{equation}{0}\setcounter{dl}{0}\setcounter{tl}{0}\setcounter{yl}{0}\setcounter{lz}{0}

 Our argument entirely relies on  Bailey's
very-well-poised$\,_6\psi_6$ summation formula of bilateral
$q$-series.
\begin{yl}[Bailey's very-well-poised$\,_6\psi_6$
summation formula]\label{tl3}(cf.\cite[II.33]{10})
\begin{align*}
\,_6\psi_6\left[%
\begin{array}{ccccccc}
 q\sqrt{a},&-q\sqrt{a}, & b, &c, & d,&e\\
 \sqrt{a},&-\sqrt{a}, & aq/b, & aq/c, & aq/d,& aq/e\\
\end{array};q,\frac{a^2q}{bcde}\right]\\
=\frac{\poq{q,aq,q/a,aq/bc,aq/bd,aq/be,aq/cd,aq/ce,aq/de}{\infty}}
{\poq{aq/b,aq/c,aq/d,aq/e,q/b,q/c,q/d,q/e, a^2q/(bcde)}{\infty}}.
 \end{align*}
 \end{yl}

 \emph{Proof of Theorem \ref{1.3}. }
First, let's define for $0<|cde|<|ab q|$ that
\begin{align}
\rho_0(a,b;c,d,e)=\sum_{k=0}^{\infty}
\left(1-\frac{aq^{2k+1}}{b}\right)
\frac{\poq{-q/b,-aq/c,-aq/d,-aq/e}{k}}
{\poq{-aq}{k}\poq{-c/b,-d/b,-e/b}{k+1}} \left(\frac{cde}{ab
q}\right)^k.\label{equa1}
  \end{align}
Then according to  (\ref{conven})
\begin{eqnarray}(a;q)_{-n}=
\frac{1}{\poq{q/a}{n}}\left(-\frac{q}{a}\right)^nq^{\binom{n}{2}}\label{nonnegative}\end{eqnarray}
we can proceed to the following computation
\begin{eqnarray*}
&&\rho_0(b,a;c,d,e)\\
&=&\sum_{k=0}^{\infty} \left(1-\frac{bq^{2k+1}}{a}\right)
\frac{\poq{-q/a}{k}\poq{-bq/c,-bq/d,-bq/e}{k}}
 {\poq{-bq}{k}\poq{-c/a,-d/a,-e/a}{k+1}}
\left(\frac{cde}{abq}\right)^k\\
&=&c_1\sum_{k=0}^{\infty}\left(1-\frac{bq^{2k+1}}{a}\right)
\frac{\poq{-q/a}{k}\poq{-b/c,-b/d,-b/e}{k+1}} {\poq{-bq}{k}
\poq{-c/a,-d/a,-e/a}{k+1}}
\left(\frac{cde}{abq}\right)^k\\
&\stackrel{\star}{=}&c_1\sum^{-1}_{k=-\infty}
 \left(1-\frac{bq^{-2k-1}}{a}\right)
\frac{\poq{-q/a}{-k-1}\poq{-b/c,-b/d,-b/e}{-k}}
{\poq{-bq}{-k-1}\poq{-c/a,-d/a,-e/a}{-k}}
\left(\frac{cde}{abq}\right)^{-k-1}\\
&=&c_1c_2\sum^{-1}_{k=-\infty}\left(1-\frac{bq^{-2k-1}}{a}\right)
\frac{\poq{-1/b}{k+1} \poq{-aq/c,-aq/d,-aq/e}{k}}
{\poq{-a}{k+1}\poq{-cq/b,-dq/b,-eq/b}{k}}
\left(\frac{cdeq}{ab}\right)^{k}\\
&=&c_1c_2c_3\sum^{-1}_{k=-\infty}\left(1-\frac{aq^{2k+1}}{b}\right)
\frac{\poq{-q/b,-aq/c,-aq/d,-aq /e}{k}}
{\poq{-aq}{k}\poq{-c/b,-d/b,-e/b}{k+1}}
\left(\frac{cde}{abq}\right)^{k},
  \end{eqnarray*}
where the line marked with ``$\star$" has been justified by the
replacement $k\to-k-1$ and the three constant $c_i$ with $i=1,2,3$
are given respectively by
\begin{eqnarray*}
  c_1 &=&\frac{1}{(1+b/c)(1+b/d)(1+b/e)},\,\, c_2= \frac{a^2b^2q}{cde}; \\
  c_3 &=&-\frac{b(1+1/b)(1+c/b)(1+d/b)(1+e/b)}{aq(1+a)}.
\end{eqnarray*}
It is easy to verify that \[c_1c_2c_3=-\frac{a(1+b)}{b(1+a)}.\]
Consequently, we  find that
 \begin{align}
&&-\frac{b(1+a)}{a(1+b)}\rho_0(b,a;c,d,e)\label{equa2}\\
&&=\sum^{-1}_{k=-\infty}\left(1-\frac{aq^{2k+1}}{b}\right)
\frac{\poq{-q/b,-aq/c,-aq/d,-aq/e}{k}}
{\poq{-aq}{k}\poq{-c/b,-d/b,-e/b}{k+1}}
\left(\frac{cde}{abq}\right)^{k}.\nonumber
   \end{align}
Observing that the right members of (\ref{equa1}) and (\ref{equa2})
form a bilateral $q$-series,  we have
\begin{eqnarray*}
&&\rho_0(a,b;c,d,e)-\frac{b(1+a)}{a(1+b)}\rho_{0}(b,a;c,d,e)
=\frac{(1-aq/b)}{(1+c/b)(1+d/b)(1+e/b)}\\
&&\times\,_6\psi_6\left[%
\begin{array}{ccccccc}
 q\sqrt{aq/b},&-q\sqrt{aq/b}, & -q/b, &-aq/c, & -aq/d,&-aq/e\\
 \sqrt{aq/b},&-\sqrt{aq/b}, & -aq, & -cq/b, & -dq/b,& -eq/b\\
\end{array};q,\frac{cde}{abq}\right]
.
  \end{eqnarray*}
Recalling Bailey's
 very-well-poised $\,_6\psi_6$ summation formula
displayed in Lemma \ref{tl3},  we  get
  \begin{eqnarray*}
&&\rho_0(a,b;c,d,e)-\frac{b(1+a)}{a(1+b)}\rho_{0}(b,a;c,d,e)
\\
&=&\frac{\poq{q,aq/b,b/a,c,d,e,cd/(ab),ce/(ab), de/(ab) }{\infty}}
{\poq{-aq,-b,-c/b,-d/b,-e/b,-c/a,-d/a,-e/a, cde/(abq) }{\infty}}.
  \end{eqnarray*}
Multiplying  both sides of this identity by $a(1+b)$ and simplifying
the resulting identity, we finally get
  \begin{align}
a(1+b)\rho_0(a,b;c,d,e)-b(1+a)\rho(b,a;c,d,e)\label{eq1}\\
=(a-b)\frac{\poq{q,aq/b,bq/a,c,d,e,cd/(ab),ce/(ab), de/(ab)
}{\infty}} {\poq{-aq,-bq,-e/a,-e/b,-c/b,-d/b,-c/a,-d/a, cde/(abq)
}{\infty}}.\nonumber
   \end{align}
 Define   \[\rho(a,b;c,d,e)=\left(1+\frac{1}{b}\right)\rho_0(a,b;c,d,e).\]
Dividing  both sides of \,(\ref{eq1}) by $ab$ and rewriting the
left-hand side of the resulting expression in terms of
$\rho(a,b;c,d,e)$, we have the formula stated in the theorem. \qed

We remark that once letting $e\mapsto 0$ in Theorem \ref{1.3}  with
 the expressions of $\rho(a,b;c,d,e)$ given in the next section,
 we get Theorem
\ref{1.2} while letting $c,d,e\mapsto 0$ simultaneously, Theorem
\ref{1.1} follows. The limiting procession is guaranteed by the
convergent condition.

\section{Applications}\setcounter{equation}{0}\setcounter{dl}{0}\setcounter{tl}{0}\setcounter{yl}{0}\setcounter{lz}{0}

By finding
 new representations for both $\rho(a,b;c,d)$ and
$\rho(a,b;c,d,e)$, we shall establish, in this section, some new
$q$-series identities. This will be realized by following Kang's
approach \cite{kang} and employing  Watson's
   $q$-analogue of Whipple's transformation formula between
  ${}_{8}W _{7}$ and ${}_{4}\phi _{3}$ series (cf.\cite[III.17]{10}).
  \begin{yl}\label{saa2} Let $a,b,c,y,z,w$ be
  such complex numbers that the following $\,_8W_7$ series is convergent
  and the $\,_4\phi_3$ series is terminating. Then it holds
  \begin{align}
   \,_8W_7\left(%
  a;  b,   c,   y, z, w
;q,\frac{a^2q^{2}}{bcyzw} \right)
=\frac{(aq,aq/(yz),aq/(yw),aq/(zw);q)_{\infty}}
{(aq/y,aq/z,aq/w, aq/(yzw);q)_{\infty}}\nonumber\\
\times\,_4\phi_3\left[%
\begin{array}{cccccccc}
  aq/(bc), & y, & z, & w \\
  & aq/b, & aq/c,& yzw/a \\
\end{array};q,q
\right].\label{Watson}
   \end{align}
   \end{yl}
Observe that  the limiting case $n\mapsto\infty$ of Watson's
transformation (\ref{Watson}), under the specification that
\[b\mapsto c/b, c\mapsto aq/c, w\mapsto q^{-n},\]
turns out to be \begin{align} &&\sum_{k=0}^{\infty}
\frac{\poq{b,y,z}{k}}{\poq{q,c,abq/c}{k}}
\left(\frac{aq}{yz}\right)^k= \frac{\poq{aq/y,aq/z}{\infty}}
 {\poq{aq,aq/(yz)}{\infty}}\label{sorry}\\
&&\times\sum_{k=0}^{\infty}(-1)^kq^{\binom{k}{2}}\frac{(1-aq^{2k})}{(1-a)}
\frac{\poq{a,c/b,aq/c,y,z}{k}} {\poq{q,abq/c,c,aq/y,aq/z}{k}}
\left(\frac{abq}{yz}\right)^k.\nonumber
  \end{align}
For  $z=q$, this becomes\begin{yl}\label{333} For
$\max\{|a/y|,|ab/y|\}<1$, there holds
 \begin{align}\label{two}
&&\sum_{k=0}^{\infty} \frac{\poq{b,y}{k}}{\poq{c,abq/c}{k}}
\left(\frac{a}{y}\right)^k= \\
&&\times\sum_{k=0}^{\infty}(-1)^kq^{\binom{k}{2}}\frac{1-aq^{2k}}{1-a/y}
\frac{\poq{c/b,aq/c,y}{k}} {\poq{abq/c,c,aq/y}{k}}
\left(\frac{ab}{y}\right)^k.\nonumber
  \end{align}
 \end{yl}
Now, we are in a position to show
\begin{dl}
 Let $\rho(a,b;c,d)$ be the same as in Theorem \ref{1.2}. Then
\begin{eqnarray}&&\rho(a,b;c,d)\nonumber\\
&&=\sum_{k=0}^{\infty}q^{\binom{k}{2}} (1-aq^{2k+1}/b)
\frac{\poq{-1/b}{k+1}\poq{-aq/c,-aq/d}{k}}
{\poq{-aq}{k}\poq{-c/b,-d/b}{k+1}} \left(\frac{cd}{b}\right)^k\label{3b}\\
&&=\left(1+\frac{1}{b}\right)
  \sum_{k=0}^{\infty}(-1)^kq^{\binom{k+1}{2}}
\frac{(1+cdq^{2k}/b)\poq{c,d,cd/(ab)}{k}}
{\poq{-aq}{k}\poq{-c/b,-d/b}{k+1}} \left(\frac{a}{b}\right)^k.
\label{3cc}
 \end{eqnarray}
\end{dl}
\pf\ Define
\[h(a,b;c,d)=\sum_{k=0}^{\infty} \frac{\poq{c,-aq/d}{k}}
 {\poq{-aq}{k}\poq{-c/b}{k+1}}
\left(-\frac{d}{b}\right)^k.\] By making the substitutions
\[
\left\{
  \begin{array}{ll}
    a &\mapsto aq/b \\
                 b &\mapsto c \\
                 c &\mapsto -aq\\
                 y &\mapsto -aq/d
  \end{array}
\right.\qquad\mbox{and}\qquad\left\{
               \begin{array}{ll}
                 a &\mapsto -cd/b \\
                 b &\mapsto -aq/d \\
                 c &\mapsto -aq\\
                 y &\mapsto c
               \end{array}
             \right.
\]
 in (\ref{two}) and then dividing
 the resulting identities  by $1+c/b$, we obtain
\begin{eqnarray*}
h(a,b;c,d)&=&\sum_{k=0}^{\infty}q^{\binom{k}{2}} (1-aq^{2k+1}/b)
\frac{\poq{-q/b,-aq/c,-aq/d}{k}}
{\poq{-aq}{k}\poq{-c/b,-d/b}{k+1}} \left(\frac{cd}{b}\right)^k\\
&=&\sum_{k=0}^{\infty}(-1)^kq^{\binom{k+1}{2}} (1+cdq^{2k}/b)
\frac{\poq{c,d,cd/(ab)}{k}} {\poq{-aq}{k}\poq{-c/b,-d/b}{k+1}}
\left(\frac{a}{b}\right)^k.
  \end{eqnarray*}
Keeping in mind of the fact that
 \[\rho(a,b;c,d)=\left(1+\frac{1}{b}\right)h(a,b;c,d),\]
we get the desired result. \qed

Obviously, \,(\ref{3b}) was missed by Kang and not recorded in
Theorem \ref{1.2}. On taking the convergent  conditions of Watson's
transformation into account, we obtain an alternative representation
for $\rho(a,b;c,d,e)$.

\begin{dl}\label{3.2}
 Let $\rho(a,b;c,d,e)$ be the same as in Theorem \ref{1.3}, at least one of the parameters $c,-aq/d,-aq/e$ be of
 the form $q^{-m},m\geq 0$. Then
\begin{align}\rho(a,b;c,d,e)=\frac{1+b}{(b+c)(1-de/(abq))}\sum_{k=0}^\infty \frac{\poq{ c,-aq/d, -aq/e}{k}}{
  \poq{-aq,  -cq/b, abq^2/(de)}{k}}q^k. \label{3ccc}
 \end{align}
\end{dl}
\pf Note that in this case,  we are able to apply Watson's
transformation to Theorem \ref{1.3}  in order to get
\begin{eqnarray*}
\rho(a,b;c,d,e)&=&c_4\,_8W_7\left(%
  aq/b; -q/b,  -aq/c,-aq/d,-aq/e,q
  ;q,\frac{cde}{abq}
\right)\\
&=&c_4c_5\sum_{k=0}^\infty \frac{\poq{ c,-aq/d, -aq/e}{k}}{
  \poq{-aq,  -cq/b, abq^2/(de)}{k}}q^k,
 \end{eqnarray*}
 where the two constant $c_4$ and $c_5$ are defined by \[c_4=\frac{(1+1/b)(1-aq/b)}{(1+c/b)(1+d/b)(1+e/b)},
\quad c_5=\frac{(1+d/b)(1+e/b)} {(1-aq/b)(1-de/abq)}\]with their
product equal to
\[c_4c_5=\frac{1+b}{(b+c)(1-de/(abq))}.\]
Therefore Theorem \ref{3.2} follows. \qed

 A few special cases of interest may be displayed as follows.
 \begin{tl} For two integers
 $r,s\geq 0$, it holds
 \begin{align}
&&\frac{1+b}{c+b}\sum_{k=0}^r\begin{bmatrix}r+s-k\\r-k\end{bmatrix}_q
\frac{\poq{ c, aq^{-s}/b}{k}}{
  \poq{-aq,-cq/b}{k}}q^{(s+1)k}\nonumber\\
  &&\qquad\quad-\frac{1+a}{c+a}\sum_{k=0}^s \begin{bmatrix}r+s-k\\s-k\end{bmatrix}_q\frac{\poq{ c,bq^{-r}/a}{k}}{
  \poq{-bq,-cq/a}{k}}q^{(r+1)k}\label{eq11}\\&&=
  \left(\frac{1}{b}-\frac{1}{a}\right)
  \frac{\poq{aq/b}{r}\poq{bq/a}{s}\poq{c}{1+r+s}} {
  \poq{-aq}{r}\poq{-bq}{s}\poq{-c/b}{r+1}\poq{-c/a}{s+1}},\nonumber  \end{align}
where the $q$-binomial coefficient
$\begin{bmatrix}n\\k\end{bmatrix}_q=\displaystyle\frac{\poq{q}{n}}{\poq{q}{k}\poq{q}{n-k}}.$
\end{tl}
\pf It follows from combination of Theorem \ref{1.3} and Theorem
\ref{3.2} with $d=-aq^{1+r}$ and $ e=-bq^{1+s}$. We are not going to
produce the tedious simplification involved.\qed

 Putting $b\mapsto \infty$ (or $b\mapsto 0$ ) in \,(\ref{eq11}), then a curious $q$-series identity follows.
 \begin{tl} Let $a\neq 0$. Then for any two  integers $r,s\geq 0$,
 there holds
 \begin{align}
\sum_{k=0}^r \begin{bmatrix}r+s-k\\r-k\end{bmatrix}_q \frac{\poq{
c}{k}}{
  \poq{-aq}{k}}q^{(s+1)k}\nonumber\\-\frac{1+a}{c+a}\sum_{k=0}^s \begin{bmatrix}r+s-k\\s-k\end{bmatrix}_q\frac{\poq{ c}{k}}{
  \poq{-cq/a}{k}}\left(-\frac{1}{a}\right)^{k}\label{eq2}\\ =
  \left(-\frac{1}{a}\right)^{s+1}
  \frac{\poq{c}{1+r+s}} {
  \poq{-aq}{r}\poq{-c/a}{s+1}}.\nonumber
   \end{align}
\end{tl}
 Further, if we let $a$
tend to zero in (\ref{eq2}), then we have
\begin{tl}Let $c,q\neq 0$.  Then for any two  integers $r,s\geq 0$, there holds
 \begin{align}
&&\sum_{k=0}^r\begin{bmatrix}r+s-k\\r-k\end{bmatrix}_q q^{(s+1)k}
\poq{ c}{k} =\poqt{c^{-1}}{1+s}
\poq{cq^{1+s}}{r}\nonumber\\
&&+c^{-1}\sum_{k=0}^s
  \begin{bmatrix}r+s-k\\s-k\end{bmatrix}_q \poqt{c^{-1}}{k
}q^{-k}.\label{eq3}
   \end{align}
\end{tl}
Next, letting $q\mapsto 1$ in \,(\ref{eq3}) and replacing $c$ by
$1/(1-x)$, we obtain a finite series transformation of interest,
which may be considered as a  supplement to the classical Pfaff
transformation for the Gauss hypergeometric function $\,_2F_1(x)$
\cite[p.68, Theorem 2.2.5]{andrews-1}.
\begin{tl}For any two  integers $r,s\geq 0,x\neq 1$,
\begin{eqnarray}
&&\sum_{k=0}^r\begin{pmatrix}r+s-k\\s\end{pmatrix}\left(\frac{x}{x-1}\right)^k
=\frac{x^{r+s+1}}{(x-1)^{r}}+(1-x)\sum_{k=0}^s
  \begin{pmatrix}r+s-k\\r\end{pmatrix}x^k,\label{idd}
  \end{eqnarray}
  where $\begin{pmatrix}n\\k\end{pmatrix}$ stands for the usual binomial
  coefficient.
\end{tl}
We remark that the special case $r=s=n$ and $x=2$ of \,(\ref{idd})
is revealed to be  \,(1.81) by Gould \cite{gould}:
\begin{eqnarray*}
&&\sum_{k=0}^n\begin{pmatrix}2n-k\\n\end{pmatrix}2^k =2^{2n}.
  \end{eqnarray*}

  The case that $c$ is of the form $q^{-m}$ in Theorem \ref{3.2} also deserves
 our consideration.
 \begin{tl} Assume the conditions in Theorem \ref{1.3}. Let $a,b\neq -1$ and
$m$
 is a nonnegative integer. Then
 \begin{align}\label{0001}
&&\frac{1+b}{1+bq^{m}}\sum_{k=0}^m \frac{\poq{ q^{-m},-aq/d,
-aq/e}{k}}{
  \poq{-aq,  -q^{1-m}/b, abq^2/(de)}{k}}q^k\nonumber\\
  &&=\frac{1+a}{1+aq^{m}}\sum_{k=0}^m\frac{\poq{q^{-m},-bq/d, -bq/e}{k}}{
  \poq{-bq,  -q^{1-m}/a, abq^2/(de)}{k}}q^k.
  \end{align}
\end{tl}
\pf It suffices to  insert $\rho(a,b;c,d,e)$ given by Theorem
\ref{3.2} and then set $c=q^{-m}$ in Theorem \ref{1.3}. After a bit
of simplification, it yields the result as claimed.\qed

Since
\[\frac{\poq{q^{-m}}{k}}{\poq{-q^{1-m}/b}{k}}=\frac{\poq{q}{m}\poq{-b}{m-k}}{\poq{q}{m-k}\poq{-b}{m}}\left(-\frac{b}{q}\right)^k,\]
there exists the limiting case  $m\mapsto\infty$ of \,(\ref{0001}),
which may be stated as
\begin{tl}\label{tlnew} For $|a|,|b|<1$, it
holds
 \begin{align} (1+b)\sum_{k=0}^\infty
\frac{\poq{-aq/d, -aq/e}{k}}{
  \poq{-aq,abq^2/(de)}{k}}(-b)^k=(1+a)\sum_{k=0}^\infty\frac{\poq{-bq/d, -bq/e}{k}}{
  \poq{-bq,abq^2/(de)}{k}}(-a)^k.\label{tlnewid}
  \end{align}
\end{tl}
In fact, Corollary \ref{tlnew} is a generalization of the symmetric
property of the Rogers-Fine function \cite[Eq.(6.3)]{fine}:
\begin{eqnarray*} (1-b)\sum_{k=0}^{\infty}
\frac{\poq{aq/d}{k}b^k}{\poq{aq}{k}} =(1-a) \sum_{k=0}^{\infty}
\frac{\poq{bq/d}{k}a^k}{\poq{bq}{k}}.
 \end{eqnarray*}

Two special cases  of this corollary are of interest.
\begin{tl}For $|d|<|q|, a\neq -q^{m}, m\leq 0$, the following hold:
\begin{eqnarray} \sum_{k=0}^\infty
\frac{\poq{-aq/d}{k}}{
  \poq{-a}{k+1}}\left(\frac{d}{q}\right)^k&=&\frac{q}{q-d};\label{par}\\
\sum_{k=0}^\infty \frac{q^{k(k-1)/2}a^{k}}{
  \poq{-a}{k+1}}&=&1.\label{par1}
 \end{eqnarray}
\end{tl}
\pf
 It suffices to show \,(\ref{par}) since (\ref{par1}) is the case $d=0$ of it. For this, write
\[f(a,d)=\sum_{k=0}^\infty
\frac{\poq{-aq/d}{k}}{
  \poq{-a}{k+1}}\left(\frac{d}{q}\right)^k.
  \]
  Actually it holds that
 \begin{eqnarray*} f(a,d)&=&\frac{1}{1+a}\lim_{e\mapsto 0}\sum_{k=0}^\infty
\frac{\poq{-aq/d, -aq/e}{k}}{
  \poq{-aq,abq^2/(de)}{k}}(-b)^k\\
  &\stackrel{(\ref{tlnewid})}{==}&\frac{1}{1+b}\lim_{e\mapsto 0}\sum_{k=0}^\infty
\frac{\poq{-bq/d, -bq/e}{k}}{
  \poq{-bq,abq^2/(de)}{k}}(-a)^k\\
  &=&f(b,d),
  \end{eqnarray*}
  which means that $f(a,d)$ is independent of the variable $a$. Hence, we have
  \[f(a,d)=f(0,d)=\frac{q}{q-d}.\] The theorem is proved.\qed

Identity (\ref{par1}) as well as its combinatorial interpretation
was also discovered by Kang \cite[Corollary
 7.4]{kang}.

\section{Some remarks on Theorem \ref{1.1}}\setcounter{equation}{0}\setcounter{dl}{0}\setcounter{tl}{0}\setcounter{yl}{0}\setcounter{lz}{0}

We end this paper by offering a new form of Ramanujan's reciprocity
theorem, which in turn leads us to a new proof for itself and for
the two-variable  generalization of
 Watson's  \textit{quintuple product identity} given by  Berndt et al. A comprehensive survey on the history and various proofs for the
latter can be found in \cite{coo}.

\begin{dl}\label{6}For $|x|<1,|q|<|a|$, it holds\begin{eqnarray}
\xi(a,x)-\frac{q}{ax}\xi(q/x,q/a)=\frac{\poq{q,ax,q/(ax)}{\infty}}{\poq{x,q/a}{\infty}},\label{gen}
 \end{eqnarray}
 where
\begin{eqnarray} \xi(a,x)=\sum_{k=0}^{\infty}\poq{a}{k}x^k=\sum_{k=0}^{\infty}
 \frac{q^{\binom{k}{2}}(-ax)^k}{\poq{x}{k+1}}.
 \end{eqnarray}
  \end{dl}
\pf Observe that the  case $c=0$ of Jackson's transformation
\cite[III.4]{10}
\begin{eqnarray}
\sum_{k=0}^{\infty}\frac{\poq{a,y}{k}}{\poq{q}{k}}x^k
=\frac{\poq{xy}{\infty}}{\poq{x}{\infty}}
\sum_{k=0}^{\infty}\frac{q^{\binom{k}{2}}(-ax)^k}
{\poq{q}{k}}\frac{\poq{y}{k}}{\poq{xy}{k}}.\label{jackson}
 \end{eqnarray}
In particular,  when $y=q$, it follows that
\begin{eqnarray} \xi(a,x) =\sum_{k=0}^{\infty}
 \frac{q^{\binom{k}{2}}(-ax)^k}{\poq{x}{k+1}}.\label{binomne9444}
 \end{eqnarray}
Hence, by  (\ref{binomne9444}) and (\ref{nonnegative}), it is easily
found that
\begin{eqnarray}
\xi(a,x)-\frac{q}{ax}\xi(q/x,q/a)=\sum_{k=-\infty}^{\infty}
\frac{(-1)^kq^{\binom{k}{2}}}{\poq{x}{k+1}}(ax)^k.
\label{binomne944}
 \end{eqnarray}
Evaluating the last sum by Ramanujan's $\,_1\psi_1$ summation
formula
\begin{eqnarray*}
\sum_{k=-\infty}^{\infty}\frac{(-1)^kq^{\binom{k}{2}}}{\poq{x}{k+1}}(ax)^k=\frac{\poq{q,ax,q/(ax)}{\infty}}
{\poq{x,q/a}{\infty}},
 \end{eqnarray*}
 we  finally get
\begin{eqnarray*}
\xi(a,x)-\frac{q}{ax}\xi(q/x,q/a)=\frac{\poq{q,ax,q/(ax)}{\infty}}{\poq{x,q/a}{\infty}},\label{recipro}
 \end{eqnarray*}
which proves the theorem.\qed

It is also worth pointing out that Ramanujan's original reciprocity
theorem, i.e., Theorem \ref{1.1}, follows from  \,(\ref{gen}) by
setting  \(x=-aq, a=-1/b\).  Note that in this case
\[\rho(a,b)=\frac{1}{b}\xi(-1/b,-aq).\]
This fact was  used by Kang. See \cite[Eq.(3.5)]{kang}. From
(\ref{binomne944}) it is clear why almost all known analytic proofs
employed Ramanujan's $\,_1\psi_1$
 summation formula.

 The next is  a new proof of Theorem
\ref{1.1}, i.e., (\ref{gen}), without invoking Ramanujan's
$\,_1\psi_1$
 summation formula.

{\sl Proof of Theorem \ref{1.1}} At first,
   replace  $a$ instead of $y$ by $q$ in (\ref{jackson}) to get
\begin{eqnarray}
\xi(y,x)=\frac{\poq{xy}{\infty}}{\poq{x}{\infty}}\lim_{t\mapsto 0}{}_{2}\phi _{1}\left[\begin{matrix}x/t,y\\
xy\end{matrix} ; q, tq\right].\label{ittt}
 \end{eqnarray}
 On the other hand, by making the substitution  \[a\mapsto x/t, b\mapsto y, c\mapsto xy, z\mapsto tq\]
in   the three-term transformation formula
 \cite[III.31]{10},
 we obtain that
\begin{eqnarray*}
{}_{2}\phi _{1}\left[\begin{matrix}x/t,y\\
xy\end{matrix} ; q,
tq\right]&-&\frac{q}{xy}\frac{\poq{x,q^2/xy,ty}{\infty}}{\poq{q/y,xy,tq/x}{\infty}}\,{}_{2}\phi
_{1}
\left[\begin{matrix}q/ty,q/x\\
q^2/xy\end{matrix} ; q,
tq\right]\\
&=&\frac{\poq{q,q/(xy)}{\infty}}{\poq{q/y,tq/x}{\infty}},\label{gennew}
 \end{eqnarray*}
from which it follows  that
\begin{eqnarray}
\lim_{t\mapsto 0}{}_{2}\phi _{1}\left[\begin{matrix}x/t,y\\
xy\end{matrix} ; q,
tq\right]&-&\frac{q}{xy}\frac{\poq{x,q^2/xy}{\infty}}{\poq{q/y,xy}{\infty}}\lim_{t\mapsto
0}{}_{2}\phi _{1}
\left[\begin{matrix}q/ty,q/x\\
q^2/xy\end{matrix} ; q,
tq\right]\nonumber\\
&=&\frac{\poq{q,q/(xy)}{\infty}}{\poq{q/y}{\infty}}.\label{gennew}
 \end{eqnarray}
Taking  (\ref{ittt}) into account and multiplying both sides of
(\ref{gennew}) by $\poq{xy}{\infty}/\poq{x}{\infty}$ we have the
formula stated in the theorem. \qed

 With Theorem \ref{6}, it will be
easier to derive the two-variable generalization of Watson's
quintuple product identity due to Berndt et al. \cite[Theorem
3.1]{ber} .\begin{tl}
\begin{align}
\sum_{k=0}^{\infty}\frac{(1-xzq^{2k})\poq{z}{k}} {\poq{x}{k+1}}
q^{k(3k-1)/2}(-1)^k\left(x^2z\right)^k\nonumber\\
-\frac{q}{zx}\sum_{k=0}^{\infty}\frac{(1-q^{2k+2}/(xz))\poq{q/x}{k}}
{\poq{q/z}{k+1}}
q^{k(3k-1)/2}(-1)^k\left(q^3/(z^2x)\right)^k\label{30}\\
=\frac{\poq{q,zx,q/(zx)}{\infty}}{\poq{x,q/z}{\infty}}.\nonumber
 \end{align}
  \end{tl}
 \pf Recall that  the limiting case  $c\mapsto 0$ while $b=c^2$
 of
(\ref{sorry}) yields
\begin{eqnarray*}
&&\sum_{k=0}^{\infty}\frac{\poq{y,z}{k}}{\poq{q}{k}}
\left(\frac{aq}{yz}\right)^k= \frac{\poq{aq/y,aq/z}{\infty}}
 {\poq{aq,aq/(yz)}{\infty}}\label{251}\\
&&\times\sum_{k=0}^{\infty}\frac{1-aq^{2k}}{1-a}
\frac{\poq{a,y,z}{k}} {\poq{q,aq/y,aq/z}{k}}
(-1)^kq^{3\binom{k}{2}}\left(\frac{a^2q^2}{yz}\right)^k
.\end{eqnarray*} Next, make the substitution $y\mapsto q, a\mapsto
xz$ in this identity so that the both sides become respectively
\begin{eqnarray*}
\xi(z,x)&=&\sum_{k=0}^{\infty}\poq{z}{k}x^k\\
&=&\sum_{k=0}^{\infty}\frac{(1-xzq^{2k})\poq{z}{k}} {\poq{x}{k+1}}
q^{k(3k-1)/2}(-1)^k\left(x^2z\right)^k.
 \end{eqnarray*}
Substituting  the last relation in \,(\ref{gen}) leads us to the
claimed result.
 \qed

In particular when we set $x=-a,z=-aq$ in \,(\ref{30}), then
Watson's celebrated quintuple product identity \cite{coo} follows
immediately
\begin{eqnarray}\sum_{k=-\infty}^{\infty}(a^2q^{2k+1}-1)a^{3k+1}q^{k(3k+1)/2}
= (q, a, q/a; q)_{\infty}(qa^2, q/a^2; q^2)_{\infty}.\label{0}
\end{eqnarray}

\bibliographystyle{amsplain}

\begin{thebibliography}{99}

\bibitem{ad}C.Adiga, N.Anitha, On a reciprocity theorem of Ramanujan,
Tamsui Oxford J.Math.Sci.\textbf{22}(2006), 9-15.

\bibitem{andrews-1} G.E.Andrews, R.Askey, and R.Roy, Special Functions,
Encyclopedia of Mathematics and Its Applications, Vol.71, Cambridge
University Press, Cambridge, UK, 1999.

\bibitem{andrews}G.E.Andrews, Ramanujan's ¡°lost¡±
notebook.I.Partial $\theta$-functions, Adv.in Math.\textbf{41}
(1981), 137-172.

\bibitem{ber0}G.E.Andrews and B.C.Berndt, Ramanujan's Lost Notebook, Part I,
Springer-Verlag, 2005; Part II, to appear.


\bibitem{ber}B.C.Berndt, S.H.Chan, B.P.Yeap,
and A.J.Yee, A reciprocity theorem for certain q-series found in
Ramanujan's lost notebook, The Ramanujan J.\textbf{13},
No.1-3(2007), 27-37.

\bibitem{coo}S.Cooper, The quintuple product identity,
Int.J.Number Theory \textbf{2}(2006), 115-161.

\bibitem{fine}N.J.Fine, Basic Hypergeometric Series and Applications,
in Mathematical Surveys and Monographs, Vol.27, Amer.Math.Soc.,
Providence, 1988.

\bibitem{10}
G.Gasper and M.Rahman, Basic Hypergeometric Series, second edition,
Cambridge University Press, Cambridge, 2004.

\bibitem{gould}H.W.Gould, Combinatorial Identities: a standardized set of tables listing 500 binomial
coefficient summations, Privately printed, Morgantown, WV, 1972.

\bibitem{liu}Z.G.Liu, Some operator identities and $q$-series transformation
formulas, Discrete Math. \textbf{265} (2003), 119-139.

\bibitem{raman}S.Ramanujan, The Lost Notebook and Other Unpublished Papers,
Narosa, New Delhi, 1988.

\bibitem{kang}S.Y.Kang, Generalizations of Ramanujan's reciprocity
theorem and their applications, J.London Math.Soc.\textbf{75}
(2)(2007), 18-34.


\end{thebibliography}

\end{document}